\documentclass[12pt,a4paper]{article}

\usepackage{amsthm}
\usepackage{amsfonts}
 \usepackage{amsmath}
\usepackage[top=20mm, bottom=20mm, left=20mm, right=20mm]{geometry}

\def\R{\mathbb{R}}

\numberwithin{equation}{section}
\theoremstyle{plain}

\newtheorem*{corollary*}{Corollary}

\def\R{\mathbb{R}}

\date{}
\begin{document}

	\title{Deterministic diffusion in dynamical systems\\ with a tiled phase space}
	
	\author{Irina L.~Nizhnik\footnote{Institute of Mathematics NAS of Ukraine, Kyiv, Ukraine, irene@imath.kiev.ua }}
	
	\maketitle

	\begin{abstract}
	The existence of normal deterministic diffusion in dynamical systems with a two-dimensional  phase space tiled by regular triangles
	(or their unions into regular hexagons) is proven.	
	\end{abstract}


\section{Introduction}

In some dynamical systems, chaotic processes may occur, for the description of which the concept of "deterministic diffusion" \, is used. With the help of simple examples
 of one-dimensional dynamical systems (DS), we explain the mechanism of appearance of deterministic diffusion and analyze
 methods used for the quantitative characterization of these processes.

 One--dimensional dynamical systems on the entire axis with discrete time are given by the recurrence relation:
\begin{equation}\label{Nizheq:1}
	x_{n+1}=f(x_{n}),
\end{equation}
where $f(x)$  is a real function given on the entire axis and  $x_0$  is a given initial value.  Equation  (\ref{Nizheq:1})
determines a trajectory
$(x_0,x_1,...,x_n,...)=x$  in the dynamical system
(\ref{Nizheq:1}) according to the initial value $x_0$ and
the form of the function $f.$
For dynamical systems admitting the chaotic behavior of trajectories and for which the function $f(x)$ has many discontinuities, the problem of
construction of the entire trajectory or even of determination of the values of $x_n$ for large $n$ is quite complicated
because, as a rule, the calculations are performed with a certain accuracy and the dependence of the subsequent
values of $x_n$ on the variations of the previous values is unstable. In addition, from the physical point of view,
the initial value $x_0$ is also speciﬁed with a certain accuracy.
Therefore, for the investigation of the behavior
of trajectories for large values of time, instead of the analysis of the evolution of system (\ref{Nizheq:1}), we can study the
evolution of measures caused by this evolution and is deﬁned on the entire straight line. If a probability measure $\mu_{0}$
(i.e., a normalized measure for which the measure of the entire axis is equal to 1) with density $\rho_{0}:$\,\,$\mu_{0}(A)=\int\limits_{A}\,\rho_{0}(x)\,dx,$
  is given at the initial time, then, for
  a unit  time, system (\ref{Nizheq:1}) maps this measure into $\mu_{1}:$\,\,$\mu_{1}(A)=\mu_0(f^{-1}(A))$, where  $f^{-1}(A)$ -- is
  the complete preimage of the set $ A$ under the mapping $f.$ The operator that maps the measure $\mu_{0}$
  into the measure $\mu_{1}$ is called the Perron–Frobenius operator. The Perron–Frobenius operator $\mathcal{F}$  
 is  linear even for nonlinear dynamical systems (\ref{Nizheq:1}) and maps the density  $\rho_{0}$ of the initial
measure into the density  $\rho_{1}$ of the measure $\mu_{1}$ as the integral operator with singular kernel
containing the Dirac delta function:
 $$\rho_{1}(x)=\mathcal{F}\rho_{0}(x)=\int\,\delta(x-f(y))
\rho_0(y)\,dy=\sum_{y_k \in f^{-1}(\{x\}) }\,
\frac{1}{|f'(y_k)|}\rho_0(y_k).$$
The investigation of the asymptotic behavior of the density
 \begin{equation}\label{Nizheq:2}
	\rho_{n}=\mathcal{F}^{n}\rho_0,\quad  \textrm{as}\quad  n\rightarrow\infty
\end{equation}
is reduced to the investigation of the
behavior of the semigroup  $\mathcal{F}^{n}$. 
There are examples of dynamical systems (\ref{Nizheq:1}) with locally expanding maps $f$  for
which the densities  $\rho_{n}$ are asymptotically Gaussian as  $n\rightarrow\infty$ independently of the choice of the density of initial
probability measure. In this case, we say that deterministic diffusion occurs in the dynamical system (\ref{Nizheq:1}).

A detailed study of dynamical systems with deterministic diffusion in the case of one--dimensional phase space was carried out in the works \cite{NN2020, N2017}, where there are references to the main results on this subject.
We present the construction of such DS with deterministic diffusion, which will be useful for studying DS in two--dimensional space that considered in section 3.

\section{Deterministic diffusion in one-dimensional dynamical systems with a tiled phase space}

Let the entire axis $\R^1$ be splitted into disjoint, semi--closed intervals $I_k=[k-\frac{1}{2},k+\frac{1}{2})$, where $k\in Z$ are integers. The function $f(x)$ of the system (\ref{Nizheq:1}) is linear and locally independent of the interval number in each interval.  We assume that  the function $f(x)=\Lambda x$ in the interval
$I_0=[-\frac{1}{2},\frac{1}{2})$, where $\Lambda=2m+1$ is an odd number, for simplicity. In each interval $I_k$ the function has the form $f(k+x)=k+\Lambda x$, where $|x|<\frac{1}{2}$. That is, the function $f(x)$ is $1$--periodic with lift. Such a dynamical system leads to the evolution of the measure as follows. If  the measure has a constant density $\rho_k$ on each interval $I_k$ at the initial moment, and the measure of the entire axis $\displaystyle \sum_{k=-\infty}^{+\infty}\,\rho_k=1$, then the Perron--Frobenius operator $\mathcal{F}$ preserves the structure of the constancy of the densities in all intervals $I_k$. However, their density changes in such a way that the measure $\mu_k=\rho_k$ of each interval $I_k$ transforms into measures of $2m+1$ intervals $I_{k-m},...,I_{0}, I_{1},...,I_{k+m}$
with densities $\displaystyle \frac{1}{\Lambda}\rho_k$, and the initial measure transforms into the sum of  specified changes of the measures from each interval.
Let $\rho_0$ be the piecewise-constant density of the initial measure in intervals, and $\rho_{1}=\mathcal{F}\rho_0$ be the density of the measure under the action of the Perron--Frobenius operator $\mathcal{F}$, then  using convolution, the explicit form of the density is
$$\rho_{1}(x)=E*\rho_{0}=\int\limits_{-\infty}^{\infty}\,E(x-y)\rho_0(y)\,dy,
$$
where
$$E(x)=\frac{1}{\Lambda}\sum_{j=-m}^{m}\,\delta(x-j),
$$
and
$\delta(x)$ -- the Dirac delta  function.  The fact that the Perron--Frobenius operator acts as a convolution makes it possible to study the evolution of measures when passing from equality (\ref{Nizheq:2}) to the Fourier transform. Since an $n$-fold convolution leads to an $n$-fold Fourier transform of the function $E(x).$ If the initial density $\rho_0$ is 1 in the interval $I_0$, then its Fourier transform is
$$\tilde{\rho_0}(\lambda)=\int\limits_{-\infty}^{\infty}\,e^{i\lambda x}\rho_0(x)\,dx=\int\limits_{-\frac{1}{2}}^{+\frac{1}{2}}\,e^{i\lambda x}\,dx=\frac{2}{\lambda}\sin\frac{\lambda}{2}.
$$
The Fourier transform of the function $E(x)$ has the form:
\begin{gather*}
	\nonumber
	\tilde{E}(\lambda)=\int\limits_{-\infty}^{\infty}\,e^{i\lambda x}E(x)\,dx=\frac{1}{\lambda}\sum_{j=-m}^{m}\,e^{i\lambda x}=\frac{1}{\Lambda}\Bigl(1+2\sum_{j=1}^{m}\,\cos \lambda j \Bigr)=\\[5mm]
	=\frac{1}{2m+1}\Bigl(2m+1-\lambda^2\sum_{j=1}^{m}\,j^2 \Bigr)+O( \lambda^4)=1-\lambda^2\frac{\Lambda^2-1}{24}+O( \lambda^4)=
	\\=e^{-\lambda^2\frac{\Lambda^2-1}{24}}+O( \lambda^4).
\end{gather*}
Therefore, passing in the formula (\ref{Nizheq:2}) to the Fourier transform, we obtain:
$$\tilde{\rho_n}(\lambda)= ( \tilde{E}(\lambda))^n \tilde{\rho_0}(\lambda) \approx
\displaystyle e^{-n\lambda^2\frac{\Lambda^2-1}{24}}+O( \lambda^4).
$$
Using the inverse Fourier transform, we have:
\begin{gather}\label{Nizheq:3}
	\nonumber \rho_n(x)=\frac{1}{2\pi}\int\limits_{-\infty}^{\infty}\,e^{-i\lambda x} \tilde{\rho_n}(\lambda)\,d\lambda \approx \frac{1}{2\pi}\int\limits_{-\infty}^{\infty}\,\displaystyle e^{-i\lambda x} e^{-n\lambda^2\frac{\Lambda^2-1}{24}}\,d\lambda =\\
	=\frac{h}{\sqrt{\pi}}e^{-h^2x^2}, \quad  \textrm{where} \quad h^2=\frac{6}{n(\Lambda^2-1)}.
\end{gather}
This is the Gaussian density. And, therefore, the considered dynamical system has a deterministic diffusion. In  \cite{NN2020}, the case of the function $f(x)=\Lambda x$ is also studied in detail,
when the numbers $\Lambda$ are even integers, which requires dividing the axis into half-intervals $I_{k,+}=[k,k+\frac{1}{2})$ and
$I_{k,-}=[k-\frac{1}{2}, k)$.  Also the case is studied when $\Lambda>2$ are non--integer numbers that densely take up the entire half--axis.

\section{Deterministic diffusion in a dynamical system with a two-dimensional phase space tiled by regular triangles}

Consider a two-dimensional space $\R^2$, tiled by disjoint triangles $\Omega_{\alpha}$ with the same side length of 2 and heights $H=3h=\sqrt{3},$
when $\displaystyle
\R^{2}=\bigcup\limits_{\alpha}\Omega_{\alpha}$, the splitting of the perimeters of the triangles is described below.
Let the entire horizontal axis $x$ in the space $\R^2$ consist of the sides of triangles, such that the origin lies in the middle of one of the sides of the triangle. Then all the centers of triangles have the following coordinates. Triangles of the first type, those whose vertices are higher than the horizontal sides, have the coordinates of the centers $(2n,h+6lh)$ or $(2n+1,-2h+6lh)$, where $l$ is an arbitrary integer. Triangles of the second type, those whose vertices are lower than their horizontal sides, have the coordinates of the centers $(2n,-h+6lh)$ or $(2n+1, 2h+6lh)$.
We  assume that triangles of the first type include triangles with their sides, except for the two lower vertices. Triangles of the second type include triangles without their sides.
Let $x=(x_1,x_2)\in\R^2$ be an arbitrary point in a two-dimensional phase space tiled by regular triangles.
We  denote by $c(x)$ the coordinate of the  triangle center to which the point $x$ belongs.

 Dynamical system of the form (\ref{Nizheq:1}) with a stretching linear mapping, a given number $\Lambda$ and a given  triangles tiling of phase space
is defined by the following equation:
\begin{equation}\label{Nizheq:*}
	x_{n+1}=c(x_{n})+\Lambda[x_{n}-c(x_{n})],\,\,x_n=(x_{n_{1}},x_{n_2}),
\end{equation}
where $c(x_{n})$ --the center of reflection of the point $x_{n},$ which is the center of the triangle to which this point belongs.
In order for dynamical system of the form (\ref{Nizheq:*}), with the parameter $\Lambda$,
to transform measures with uniform density in a triangle tiling of the phase space into a similar one, it is necessary and sufficient that $\Lambda=4+3m $, where $m$ is an integer.
Then, the measure $\mu$ in a triangle with the center $c$ and uniform density is transformed by the dynamical system into the same measure $\mu$ with uniform density in a triangle with the same center
$c,$ but with the length of the sides in $\Lambda$ times greater. To describe such an action of  dynamical system, it is useful to divide the uniform density $\rho$
in the entire triangle tiled space  into the sum of densities $\rho^{(1)}$ and $\rho^{(2)}$  concentrated, respectively, in triangles of the first and second types $\rho=\rho^{(1)}+\rho^{(2)}$. Then, the Perron---Frobenius operator $\mathcal{F}$ transforms the initial measure $\mu_0$ with density $\rho_0$ into:
\begin{equation}\label{Nizheq:5}
	\begin{pmatrix}
		\rho_1^{(1)}\\
		\rho_1^{(2)}
	\end{pmatrix}=E *  \begin{pmatrix}
		\rho_0^{(1)}\\
		\rho_0^{(2)}
	\end{pmatrix},\quad \rho_j=\rho_j^{(1)}+\rho_j^{(2)},\,\,j=0,1,
\end{equation}
where the matrix
$E=\begin{pmatrix}
	E_{11} & E_{12} \\
	E_{21} & E_{22}
\end{pmatrix}.$
For the simplest case $\Lambda=4$, the elements of the matrix $E$
have the form:
\begin{gather*}
	E_{11}(x)=\frac{1}{16}\Bigl([\delta(x_1+3)+\delta(x_1+1)+\delta(x_1-1)+\delta(x_1-3)]\delta(x_2+3h)+ \\
	+[\delta(x_1+2)+\delta(x_1)+\delta(x_1-2)]\delta(x_2)+ \\
	+[\delta(x_1+1)+\delta(x_1-1)]\delta(x_2-3h)+\delta(x_1) \delta(x_2-6h)\Bigr),
\end{gather*}
\begin{gather*}
	E_{21}(x)=\frac{1}{16}\Bigl([\delta(x_1+2)+\delta(x_1)+\delta(x_1-2)]\delta(x_2+2h)+ \\
	+[\delta(x_1+1)+\delta(x_1-1)]\delta(x_2-h)+\delta(x_1) \delta(x_2-4h)\Bigr),
\end{gather*}
Since the Dirac delta function is an even function, then $E_{12}(x_1, x_2)=E_{21}(x_1,-x_2),$\qquad$E_{22}(x_1, x_2)=E_{11}(x_1,-x_2).$
The Fourier transform (\ref{Nizheq:5}) gives:
\begin{equation}\label{Nizheq:6}
	\begin{pmatrix}
		\tilde{\rho}_1^{(1)}\\
		\tilde{\rho}_1^{(2)}
	\end{pmatrix}= \tilde{E}  \begin{pmatrix}
		\tilde{\rho}_0^{(1)}\\
		\tilde{\rho}_0^{(2)}
	\end{pmatrix}
	\quad  \textrm{and}\quad
	\begin{pmatrix}
		\tilde{\rho}_n^{(1)}\\
		\tilde{\rho}_n^{(2)}
	\end{pmatrix}=(\tilde{E})^n  \begin{pmatrix}
		\tilde{\rho}_0^{(1)}\\
		\tilde{\rho}_0^{(2)}
	\end{pmatrix}.
\end{equation}
The Fourier transform of functions
\begin{equation}\label{Nizheq:7}
	\tilde{E}_{kj}(\lambda)=\int\limits_{\R^2}e^{i\lambda x}E_{kj}(x)\,dx, \quad  \textrm{where}\quad \lambda=(\lambda_1,\lambda_2),\,\, x=(x_1,x_2),\,\,k,j=1,2
\end{equation}
have the following principal values:
$\tilde{E}_{11}(\lambda),$\,$\tilde{E}_{22}(\lambda)$ having the same  principal  parts:
$\displaystyle \frac{5}{8}-\frac{15}{16}\lambda^2+O(\lambda^4),$ and $\displaystyle \tilde{E}_{12}(\lambda)$,\, $\tilde{E}_{21}(\lambda)$ have the form:
$\displaystyle \frac{3}{8}-\frac{5}{16}\lambda^2+O(\lambda^4).$

Therefore, the principal  part of the matrix $ \tilde{E}(\lambda)$ has the form:
$$\tilde{E}(\lambda)=
\begin{pmatrix}
\displaystyle	\frac{5}{8}-\frac{15}{16}\lambda^2\,	  & \,\displaystyle \frac{3}{8}-\frac{5}{16}\lambda^2 \\[5mm]
\displaystyle	\frac{3}{8}-\frac{5}{16}\lambda^2\, & \,\displaystyle \frac{5}{8}-\frac{15}{16}\lambda^2
\end{pmatrix}$$
and is symmetric. If we consider the matrix $$A=U \tilde{E}U^{-1}=
\begin{pmatrix}
\displaystyle	1-\frac{5}{4}\lambda^2	  &  0\\[5mm]
	0 & \displaystyle \frac{1}{4}-\frac{5}{8}\lambda^2
\end{pmatrix},$$
where $\displaystyle U=\displaystyle \frac{1}{\sqrt{2}} \begin{pmatrix}
	1	  &  1\\
	-1 &  1
\end{pmatrix}$ -- a unitary matrix, $U^{-1}=U^*.$
Therefore, when $n\rightarrow\infty$
$$A^n=\begin{pmatrix}
\displaystyle 	\Bigl(1-\frac{5}{4}\lambda^2\Bigr)^n	  &  0\\
	0 & \displaystyle \Bigl (\frac{1}{4}-\frac{5}{8}\lambda^2\Bigr)^n
\end{pmatrix}=\begin{pmatrix}
	e^{-\frac{5}{4}n\lambda^2}	  &  0\\[5mm]
	0 & 0
\end{pmatrix}+O\Bigl(\frac{1}{n}\Bigr).
$$
Then, the equality (\ref {Nizheq:6}) can be represented as folowing:
$$ U\begin{pmatrix}
	\tilde{\rho}_n^{(1)}\\
	\tilde{\rho}_n^{(2)}
\end{pmatrix}=A^n U \begin{pmatrix}
	\tilde{\rho}_0^{(1)}\\
	\tilde{\rho}_0^{(2)}
\end{pmatrix}.
$$
From these results it follows that
\begin{equation}\label{Nizheq:8}
	\tilde{\rho}_n(\lambda)=\tilde{\rho}_n^{(1)}+\tilde{\rho}_n^{(2)}=e^{-\frac{5}{4}n\lambda^2}\tilde{\rho}_0.
\end{equation}
This equality (\ref{Nizheq:8}) means that the density
$$\rho_n(x)=\frac{1}{4\pi^2}\int\limits_{\R^2}\,e^{-i\lambda x}e^{-\frac{5}{4}n\lambda^2}\,d\lambda=\frac{1}{5\pi n}e^{-\frac{x^2}{5n}}
$$
is the Gaussian two-dimensional measure. That is, in the dynamical system under  consideration, normal deterministic diffusion occurs.

If the two-dimensional space $\R^2$ is tiled with regular hexagons, each divided into six regular triangles, then the dynamical system of  (\ref{Nizheq:*}), where $c$ are the coordinates of the  hexagon centers, and $\Lambda \geq 3$ is an odd number, the initial measure, which has a constant density in a certain number of triangles, is transformed by the Perron---Frobenius operator $\mathcal{F}$ into a measure of the same structure. Applying a similar method, as when considering triangular tiling, it is proved that in the dynamical system (\ref{Nizheq:*}) with hexagon tiling, for odd numbers $\Lambda\geq 3$, deterministic diffusion occurs.

{\bf Acknowledgments.}\\
The author is  grateful  to the Simons Foundation  for the financial support of the Institute of Mathematics of the National Academy of Sciences of Ukraine
(1290607, I.N.).

\end{document}